\documentclass[letterpaper,11pt]{article}

\usepackage{amsthm}
\usepackage{amsmath}

\newtheorem{thm}{Theorem}[section]

\newtheorem{lemma}[thm]{Lemma}
\newtheorem{prop}[thm]{Proposition}

\newcommand{\beq}[1]{\begin{equation}\label{#1}}
\newcommand{\enq}[0]{\end{equation}}

\newcommand{\bn}[0]{\bigskip\noindent}
\newcommand{\mn}[0]{\medskip\noindent}
\newcommand{\nin}[0]{\noindent}

\newcommand{\sub}[0]{\subseteq}

\renewcommand{\dots}[0]{,\ldots,}

\newcommand{\A}[0]{{\cal A}}

\newcommand{\T}[0]{{\cal T}}

\newcommand{\GG}[0]{{\bf G}}
\newcommand{\HH}[0]{{\bf H}}

\newcommand{\TT}[0]{{\bf T}}

\newcommand{\E}[0]{{\sf E}}

\renewcommand{\qed}[0]{\begin{flushright} \rule{2mm}{3mm} \end{flushright}}

\newcommand{\C}[2]{{{#1}\choose{{#2}}}}
\newcommand{\Cc}[0]{\tbinom}
\newcommand{\ga}[0]{\alpha }
\newcommand{\gb}[0]{\beta }
\newcommand{\gc}[0]{\gamma }
\newcommand{\gd}[0]{\delta }

\newcommand{\gl}[0]{\lambda }

\newcommand{\gO}[0]{\Omega}

\newcommand{\gz}[0]{\zeta}
\newcommand{\eps}[0]{\varepsilon }

\newcommand{\nine}[0]{10}

\newcommand{\sugg}[1]{}

\begin{document}

\renewcommand{\thefootnote}{\fnsymbol{footnote}}
\title{Upper tails for triangles}

\author{
B. DeMarco\footnotemark $~$ and J. Kahn\footnotemark
}
\date{}
\footnotetext{ * supported by the U.S.
Department of Homeland Security under Grant Award Number 2007-ST-104-000006.}
\footnotetext{ $\dag$ Supported by NSF grant DMS0701175.}

\sugg{
\author{
B. Demarco
\footnotetext{ * supported by the U.S.
Department of Homeland Security under Grant Award Number 2007-ST-104-000006.}
and J. Kahn
\footnotetext{ $\dag$ Supported by NSF grant DMS0701175.}
}
\date{}
}

\maketitle

\begin{abstract}
With $\xi$ the number of triangles in the usual (Erd\H{o}s-R\'enyi) random
graph $G(m,p)$, $p>1/m$ and $\eta>0$, we show (for some $C_{\eta}>0$)
$$\Pr(\xi> (1+\eta)\E \xi) <
\exp[-C_{\eta}\min\{ m^2p^2\log(1/p),m^3p^3\}].$$
This is tight up to the value of $C_{\eta}$.

\end{abstract}

\section{Introduction}\label{Intro}

Write $\GG=G(m,p)$ for the usual (Erd\H{o}s-R\'enyi) random graph,
and $\xi=\xi(\GG)$ for the number of triangles in $\GG$.
The main purpose of this note is to prove

\begin{thm}\label{T1}
For any $\eta>0$ and $p> m^{-1} \ln m$,
\beq{bd1}
\Pr(\xi > (1+\eta)\Cc{m}{3}p^3) <p^{\Omega_{\eta}( m^2p^2)}.
\enq
\end{thm}

\nin
As perhaps first observed by Vu \cite{Vu1},
the bound is tight up to the value of the constant in the exponent,
since the probability that $\GG$ contains a complete graph on (say)
$2mp$ vertices is $\exp_p[O(m^2p^2)]$.
The lower bound on $p$ is needed because below this the behavior
changes, and a better lower bound on the probability in \eqref{bd1}
(ignoring the constants in the exponents)
is the more natural one given by

\begin{prop}\label{Plb}
For $p>1/m$,
$
~\Pr(\xi > 2\Cc{m}{3}p^3) >\exp[-O( m^3p^3)],
$
unless the probability in question is zero
(i.e. unless $2\Cc{m}{3}p^3\geq\C{m}{3})$.
\end{prop}

\nin
For $p<1/m$ the question is uninteresting:  one easily has
(for $\eta \in (0,1)$, say)
$\Pr(\xi >(1+\eta)\C{m}{3}p^3) = \Pr(\xi \geq 1)=\Theta(m^3p^3)$,
with the lower bound given, e.g., by (one of) the Bonferroni Inequalities.
Of course Proposition \ref{Plb} is not news if $p> m^{-1}\ln m$,
but
see the discussion in
Section \ref{Lower}.

The proof of Theorem \ref{T1} extends without serious modification to give
a full resolution of the question for all $p$ of interest, saying
that in all cases the larger
of the preceding two lower bounds is the truth:

\begin{thm}\label{T3}
For any $\eta>0$ and $p$ 
$$
\Pr(\xi > (1+\eta)\Cc{m}{3}p^3) <
\exp[-\Omega_{\eta}(\min\{ m^2p^2\log(1/p),m^3p^3\})].
$$
\end{thm}

\nin
(Again, the statement is not interesting for $p<1/m$.)

\medskip
A result close to Theorem \ref{T1}
was proved independently by S. Chatterjee
in \cite{Chat}.
The approaches there and here are related, though the proof given
here is quite a bit shorter and proves a little more, in that
\cite{Chat} requires $p>C_{\eta}^{-1} m^{-1}\ln m $ for
some small $C_{\eta}$.
(Prompted by \cite{Chat},
a version of the present paper proving only Theorem \ref{T1}---the
proof is essentially the same as that given here---was posted at
\cite{DK}.)

Though only a first case of the analogous question for copies of
a general (fixed) graph $H$,
the problem addressed by Theorem \ref{T1} (or Theorem \ref{T3})
has a surprisingly
substantial history; see e.g. \cite{Vu1,Jan4,Jan5, Kim2,JOR,Chat2}
for some of this, or
\cite{Chat} for a thorough account.
Here we just mention that Kim and Vu \cite{Kim2} were the first to prove
that the probability in question is $\exp[-\Omega_{\eta}(m^2p^2)]$.
(See Section \ref{Remarks} for a little more on general $H$.)

\medskip
For the proof of Theorem \ref{T3}
it's convenient to work with a tripartite
version.
Let $\HH$
be the random tripartite graph on $n+n+n$
vertices; thus $V=V(\HH)$ is the disjoint union of three $n$-sets,
say $V_1,V_2,V_3$, and
$\Pr(xy\in E(\HH)=p)$ whenever $x,y$ are in distinct $V_i$'s, these
choices made independently.
Then with $\xi'$ the number of triangles in $\HH$ (so $\E\xi' = n^3p^3$),
we show

\begin{thm}\label{T2}
For any $\gd>0$ and $p$, 
\beq{T4disp}
\Pr(\xi'>(1+\gd)n^3p^3)<
\exp[-\Omega_{\gd}(\min\{ n^2p^2\log(1/p),n^3p^3\})].
\enq
\end{thm}

That Theorem \ref{T2} implies Theorem \ref{T3} is presumably well-known,
but we give the easy argument.
It is of course enough to prove Theorem \ref{T3} when $m=3n$.
Let $\eta$ be as in Theorem \ref{T3} and set $\gd = \eta/(2+\eta)$.
We may choose $\HH$ by choosing $\GG$ (on $V=[3n]$)
and a uniform equipartition
$V_1\cup V_2\cup V_3$ of $V$, and setting
$$E(\HH)=\{xy\in E(\GG): x,y ~\mbox{belong to distinct $V_i$'s}\}.$$
Of course
\beq{Exi'}
\E[\xi'|G]=\rho\xi(G),
\enq
where $\rho
= n^3/\C{3n}{3}$ ($\sim 2/9$) and
the conditioning event is $\{\GG=G\}$.
On the other hand, with $\ga(G) = \Pr(\xi'< (1-\gd)\rho\xi(G)|G),$
we have
$$
\E[\xi'|G] \leq \ga(G)(1-\gd)\rho\xi(G) +(1-\ga(G))\xi(G),
$$
whence, using \eqref{Exi'}, $\ga(G)\leq 1-\rho\gd/(1-\rho+\rho\gd) =:1-\gb$.
Thus (by Theorem \ref{T2})
\begin{eqnarray*}
\exp[-\Omega_{\gd}(\min\{ n^2p^2\log(1/p),n^3p^3\})]
&>&\Pr(\xi'>(1+\gd)n^3p^3) \\
&\geq & \gb
\Pr(\xi >\tfrac{1+\gd}{1-\gd}\Cc{3n}{3}p^3),
\end{eqnarray*}
and (noting $(1+\gd)/(1-\gd) =1+\eta$) Theorem \ref{T3} follows.\qed

Theorem \ref{T2} is proved in Section \ref{Proof}.
There is more than one way to prove it, but we just give
the easiest (that we know).
Section \ref{Lower} gives the easy (and probably known) proof of
Proposition \ref{Plb}, and in
Section \ref{Remarks} we briefly mention the situation  for general $H$.
As noted there, the present method will surely prove more than
Theorem \ref{T3}, but at this writing we don't know exactly how much more.

\section{Proof of Theorem \ref{T2}}\label{Proof}

We rename the parts of our tripartition $A,B,C$ and always take $a,b,c$ to
be elements of $A,B,C$ respectively.
A triangle is then simply denoted $abc$.
The set of triangles of $\HH$ is denoted $\TT$.
As usual $N_Y(x)=\{y\in Y:xy\in E(\HH)\}$ and $d_Y(x) = |N_Y(x)|$,
and we also use (e.g.) $d(a)=\max\{d_B(a),d_C(a)\}$ and
$d(a,b)=|N_C(a)\cap N_C(b)|$.
\sugg{
We use $d(a)=\max\{d_B(a),d_C(a)\}$, where $d_Y(x)=|\{y\in Y:xy\in E(\HH)\}|$,
(and similarly for $d(b),d(c)$, and
$d(a,b)=|\{c:ac,bc\in E(\HH)\}|$.
{\bf OR}
We use $B(a)=\{b\in B:ab\in E(\HH)\}$,
$d(a)=\max\{|B(a)|,|C(a)|\}$
and
$d(a,b)=|C(a)\cap C(b)|$, and similarly with the roles of $a,b,c$
permuted.
}
For disjoint $X,Y\sub V$ we use $\nabla(X)$ (resp. $\nabla(X,Y)$)
for the set of edges with one end in $X$
(resp. one end in each of $X,Y$).

Set $t=\ln (1/p)$, $s=\min\{t,np\}$,
$\ga = \gd/3$, $\eps = .02 \ga$ and
(say)
$\gc = 1/e$.
We may assume
\beq{psmall}
p < (1+\gd)^{-1/3},
\enq
since otherwise the left side of \eqref{T4disp} is zero.
We may also assume:  $\gd$---so also $\eps$---is (fixed but) small
(since \eqref{T4disp} becomes weaker as $\gd$ grows);
given $\gd$, $n$ is large (formally, $n>n_\gd$);
and, say,
\beq{plarge}
p>\eps^{-3}n^{-1}
\enq
(since for smaller $p$,
Theorem \ref{T2} becomes trivial with an appropriate
$\Omega_{\gd}$).

\sugg{By \cite{Kim2} we may assume $p$ is not too close to 1.
(This is unnecessary but saves a little fussing with the parameters.)
We may also assume $\eps$ (or $\gd$) is (fixed but) somewhat small and $n$ is large,
and, say, $p>\eps^{-3}n^{-1}$
(since for smaller $p$
Theorem \ref{T2} becomes trivial with an appropriate choice of $\Omega_{\gd}$).}

We say an event occurs {\em with large probability} (w.l.p.)
if its probability is at least $1-\exp_p[T\eps^4 n^2p^2]$ for some fixed $T>0$
and small enough $\eps$ (and $p$ satisfying \eqref{plarge}),
and write ``$\ga<^*\gb$"
for ``w.l.p. $\ga < \gb$."
Note that an
intersection of $O(n)$ events that hold w.l.p.
also holds w.l.p.

Let $A'=\{a:d_C(a)\leq np^{1-\gc}\}$ and $B'=\{b:d_C(b)\leq np^{1-\gc}\}$.
The next three assertions imply Theorem \ref{T2}.
\beq{Cl1}
\mbox{w.l.p. $|\{abc\in \TT:a\not\in A'$ or $b\not\in B'$}\}| < \ga n^3p^3;
\enq
\beq{Cl2}
\mbox{w.l.p. $|\{abc\in \TT:a\in A',b\in B',
d(a,b) > \nine np/s\}| < \ga n^3p^3$};
\enq
\beq{Cl3}
\Pr(|\{abc\in \TT:
d(a,b) \leq \nine np/s\}| > (1+\ga) n^3p^3)
< \exp[-\Omega_{\gd}(n^2p^2s)].
\enq

\medskip
We use $B(m,\ga)$ for a random variable
with the binomial distribution ${\rm Bin}(m,\ga)$.
For the proofs of \eqref{Cl1}-\eqref{Cl3} we need
the following standard bounds;
see e.g. \cite[Theorem A.1.12]{Alon2}, \cite[Theorem 2.1(a)]{JLR}
and \cite[Lemma 8.2]{Beck-Chen}.

\begin{lemma}\label{L1}
There is a fixed $C>0$ so that for any $K>1+\eps$, m and $\ga$,
\beq{bin}
\Pr(B(m,\ga) \geq Km\ga) <\min\{(e/K)^{Km\ga}, \exp[-C\eps^2Km\ga]\}.
\enq
\end{lemma}

\nin
(Though unnecessary, the $K$ in the second expression
will be helpful below.)
When $m=n$ and $\ga=p$ we use $q_K$ for the r.h.s. of \eqref{bin}.

\begin{lemma}\label{L2}
Suppose $w_1\dots w_m \in [0,z]$.
Let $\xi_1\dots \xi_m$ be independent Bernoullis,
$\xi = \sum \xi_i w_i$, and $\E\xi =\mu$.
Then for any $\eta >0$ and $\gl\geq \eta \mu$,
$$\Pr(\gz > \mu+\gl) < \exp [- \gO_\eta(\gl/z)].$$
\end{lemma}

\begin{lemma}\label{Cor}
For $K>1+\eps$ and $X\in \{A,B,C\}$,
\beq{cor}
|\{x\in X:d(x)\geq Knp\}|<^* r_K:=\left\{\begin{array}{ll}
3n\cdot \eps K^{-4} &\mbox{if $q_K>n^{-2}$}\\
\tfrac{\eps^2 npt}{K\ln K} &\mbox{otherwise}.
\end{array}\right.
\enq
\end{lemma}

\nin
The first, {\em ad hoc} value of $r_K$ is for use
in the proof of \eqref{Cl3}.
Note that
\beq{corbds}\frac{\eps^2npt}{K\ln K} < \left\{\begin{array}{ll}
2\eps npt/K&\mbox{if $K > 1+\eps$}\\
\eps np/K&\mbox{if $K>p^{-\eps}$.}
\end{array}\right.
\enq

\mn
{\em Proof of Corollary} \ref{Cor}.
Write $N$ for the left side of \eqref{corbds}
and let $q_K=q, r_K=r$ and (w.l.o.g) $X=A$.  If $q\leq n^{-2}$ then,
since the $d_B(a)$'s and $d_C(a)$'s are independent copies of $B(n,p)$,
two applications of Lemma \ref{L1} (and a little checking) give
$$\Pr(N\geq r)<\Pr(B(2n,q)\geq\lceil r\rceil)
<(2e\sqrt{q})^r <\exp[-\gO(\eps^4n^2p^2t)].$$

If $q>n^{-2}$ then $\exp[-C\eps^2Knp]>q$
implies $Knp<2C^{-1}\eps^{-2}\ln n$, while \eqref{plarge} gives
$q\leq \exp[-C\eps^2Knp]< \exp[-C\eps^{-1}K]<\eps K^{-4}$
(the last inequality gotten by observing that
$\exp[C\eps^{-1}K]\eps K^{-4}$
is minimized at $K=4\eps/C$ and assuming, as we may,
that $\eps<(Ce/4)^{4/3}$).  It follows that
$$\Pr(N\geq r)<\Pr(B(2n,q)\geq r)<\exp[-\Omega(\eps n K^{-4})]<p^{\Omega(n^2p^2)},$$
where the second inequality uses $r>3nq/2$ (and Lemma \ref{L1})
and the (very crude) third inequality uses the above upper bound on $Knp$.
\qed

\nin
We will also make occasional use of the fact that for any $\gb>0$ and $p$,
\beq{max}
p^\gb \ln(1/p) \leq (e\gb)^{-1} ~~~\mbox{and}~~~
p^\gb \ln^2(1/p) \leq 4(e\gb)^{-2}.
\enq

\mn
{\em Proof of} \eqref{Cl1}.
For $K>p^{-\gc}$
($> 1+\eps$; see \eqref{psmall})
Lemma \ref{Cor} (with \eqref{corbds}) gives
$
|\{a:d(a)>Knp\}|<^* \eps np/K
$
(note $K>p^{-\gc}$ implies $q_K<n^{-2}$),
and similarly with $b$ in place of $a$.
On the other hand,
with $K_a=d(a)/(np)$,
\beq{2.2}
\mbox{w.l.p.} ~~~|\nabla(N_B(a),N_C(a))| <
\max\{2K_a^2n^2p^3, n^2p^2t\} =:\gb_a ~~~\forall a
\enq
(and similarly for $b$), since, given any $\nabla(A)$,
the probability that the event in
\eqref{2.2} fails is (again using Lemma \ref{L1})
less than
$$
\sum_a\Pr(B(K_a^2n^2p^2,p)>\gb_a) <
\sum_a\exp[-\Omega( \gb_a)] < \exp[-\Omega(n^2p^2t)].
$$
This gives \eqref{Cl1} since,
with $J=\sqrt{t/(2p)}$ (so $\gb_a= 2K_a^2n^2p^3$ iff $K_a\geq J$)
and $u=\lfloor-\log_2 (Jp)\rfloor$,
$|\{abc\in\TT:a\not\in A'\}|$
is at most
$$
|\{abc\in \TT:p^{-\gc}\leq K_a\leq J\}|
+ \sum_{i=0}^u |\{abc\in \TT:K_a\in [2^i J,2^{i+1}J]\}|
~~~~~~
$$
$$
~~~~~~~~~~~~~~~~~~~~~~~~
<^* \eps np^{1+\gc}n^2p^2t + \sum_{i=0}^u\tfrac{\eps np}{2^iJ}
\cdot 2\cdot 2^{2i+2}J^2n^2p^3 ~<~ 17\eps n^3p^3
$$
(using \eqref{max}), and, of course, similarly for
$|\{abc\in\TT:b\not\in B'\}|$.\qed

\mn
{\em Proof of} \eqref{Cl2}.
For $K\geq J:=\nine /s$, let
$A_K =\{a:\exists b\in B', ~d(a,b)>Knp\}$, and define $B_K$ similarly.
Given $\nabla(B,C)$ the events $\{a\in A_K\}$ are independent with,
for each $a$,
$$
\Pr(a\in A_K)~ <~ n\Pr(B(np^{1-\gc},p)>Knp) ~<~ np^{Knp/2}
=:q,
$$
using Lemma \ref{L1} (with
$ep^{1-\gc}/K<p^{1/2}$,
which follows from \eqref{max})
for the second inequality.
Now
$Knpt \geq 10\max\{np,t\}> 7\ln n$ (say) implies both
$enq^{1/2}<1$ and $q< p^{Knp/4}$, so we have
(again using Lemma \ref{L1})
$$\Pr(|A_K|\geq \eps np/K)
< (enq/\lceil \eps np/K\rceil)^{\eps np/K}
< (q^{1/2})^{\eps np/K}<\exp [-\eps n^2p^2t/8].$$
Thus $|A_K|<^* \eps np/K$,
and similarly for $B_K$.

\sugg{Since $1>enq^{1/2}$ (see \eqref{plarge}),
$$\Pr(|A_K|\geq \eps np/K)
< (enq/\lceil \eps np/K\rceil)^{\eps np/K}
< (q^{1/2})^{\eps np/K}<\exp [-\eps n^2p^2t/8]$$
(using Lemma \ref{L1} for the first inequality
and $Knpt \geq 10\max\{np,t\}> 4\ln n$ for the third).
Thus
$|A_K|<^* \eps np/K$,
and similarly for $B_K$.
}

Now thinking of first choosing
$\nabla(C)$ (which determines the $A_K$'s and $B_K$'s), we have
$|A_J|,|B_J|<^* \eps nps$, so that
$
\E |\nabla(A_J,B_J)| <^*\eps^2s^2n^2p^3.
$
Lemma \ref{L1}
(using, say, $\eps s^2p<p^{1/2}$, which follows from \eqref{max}),
then gives
$$
|\nabla(A_J,B_J)| <^* \eps n^2p^2.
$$
We may then bound the left side of \eqref{Cl2} by
$$
|\nabla(A_J,B_J)|np + \sum_{i\geq 0} 2^{1-i}\eps^2n^3p^3
<^* (\eps +4\eps^2)n^3p^3,
$$
where the first term corresponds to $abc$'s with
$d(a,b)\in [Jnp,np]$, and the $i$th summand to
those with $d(a,b)\in [2^i np,2^{i+1}np]$
(using $|\{abc\in \TT:a\in A',b\in B', d(a,b)\in[Knp,2Knp]\}|\leq
|A_K||B_K|2Knp <^* 2\eps^2 n^3p^3/K$).\qed

\mn
{\em Proof of} \eqref{Cl3}.
We first show
\beq{d2c}
\sum \{d^2(c):d(c)>(1+\eps)np\} <^* 40\eps n^3p^2.
\enq
Setting $v=\lfloor -\log_2 ((1+\eps)p)\rfloor$,
$u = \lfloor -\log_2((1+\eps)p^\eps)\rfloor$,
and using Lemma \ref{Cor} (with \eqref{corbds}),
we have
$$
\mbox{$\sum \{d^2(c):\tfrac{d(c)}{(1+\eps)np}\in [2^i,2^{i+1}]\}$} <^*
\left\{\begin{array}{ll}
4\eps(1+\eps) n^3p^32^i&\mbox{if $i> u$}\\
8\eps(1+\eps)n^3p^3t2^i&\mbox{if $u\geq i\geq 0$,}
\end{array}\right.
$$
provided $K(i):=(1+\eps)2^i$ satisfies
$q_{K(i)} \leq n^{-2}$.
The left side of \eqref{d2c} is thus w.l.p. at most
$$
3\eps n^3p^2\sum_{i\geq 0}2^{-2i+2}
+ 4\eps(1+\eps) n^3p^3 [2t\sum_{i=0}^u 2^i +\sum_{i=u+1}^v 2^i]
< 40 \eps n^3p^2,
$$
where the first term on the left, covering $c$'s with
$\tfrac{d(c)}{(1+\eps)np}\in [2^i,2^{i+1}]$ for an $i$ with
$q_{K(i)} > n^{-2}$,
again comes from Lemma \ref{Cor},
and we used \eqref{max} to say (say)
$p^{3-\eps}t<p^2$.

\medskip
Finally, set $\xi_{ab}={\bf 1}_{\{ab\in E(\HH)\}}$.
We have
\beq{sumdab}
\sum\{d(a,b):d(a,b)\leq \nine np/ s\}\leq
\sum d^2(c)<^* (1+42\eps+\eps^2)n^3p^2
\enq
(by \eqref{d2c}, where ``$<^*$" refers to the choice of $\nabla(C)$),
and
$$
|\{abc\in \TT:
d(a,b) \leq \nine np/s\}| =
\sum\{\xi_{ab}d(a,b):d(a,b)\leq \nine np/ s\},
$$
so that
Lemma \ref{L2} (with $z= \nine np/s$) combined with \eqref{sumdab} gives
$$
\Pr(|\{abc\in \TT:
d(a,b) \leq \nine np/s\}| > (1+\ga)n^3p^3) < \exp[-\Omega_{\ga}( n^2p^2s)].
$$\qed

\section{Lower bound for small $p$}\label{Lower}

Here we prove Proposition \ref{Plb}.
As noted in Section \ref{Intro}, we only need to prove it for
$p< m^{-1}\ln m$;
but
we will show a little more:

\begin{prop}\label{Plb'}
For $1/m<p<o(m^{-5/6})$,
$$\Pr(
\mbox{$\GG$ contains
$2\Cc{m}{3}p^3$ disjoint triangles}\}
>\exp[-O(m^3p^3)].$$
\end{prop}

\nin
The point is that the natural
bound of Proposition \ref{Plb} is achieved,
at least for $p$ as in Proposition \ref{Plb'},
more ``generically"
than the usually stronger $\exp_p[O(m^2p^2)]$.
There are similar statements for larger $p$---for example,
for $1/m<p<o(m^{-1/2})$
Proposition \ref{Plb'} is true with ``disjoint"
replaced by ``edge-disjoint"---but we will not pursue this tangent
further here.

\mn
{\em Proof.}
Write $\T=\T_m$ for the set of triangles of $K_m$ and set
$M=\lceil 2\C{m}{3}p^3\rceil$.
The number of sets $S$ consisting of $M$ vertex-disjoint members of $\T$ is
\beq{Scount}
\tfrac{1}{M!}\prod_{i=0}^{M-1}\Cc{m-3i}{3}\sim
\tfrac{1}{M!}\Cc{m}{3}^M
\enq
(where ``$\sim$" uses $p <o(m^{-5/6})$).
For such $S$ let $Q_S$ and $R_S$ be the events
$\{\mbox{$\GG$ contains all triangles of $S$}\}$
and
$\{\mbox{$S$ is the set of triangles of $\GG$}\}$;
let $a=3M(m-3)$ be the number of members of $\T$ that share
a (necessarily unique) edge with the union of the triangles of $S$,
and set $b= \C{m}{3}-M-a$.
Then
\beq{PrRS}
Pr(R_S) = \Pr(Q_S)\Pr(R_S|Q_S)
\geq  p^{3M} (1-p^2)^a (1-p^3)^b,
\enq
where the inequality is given by Harris' Inequality \cite{Harris}
(which for our purposes says that for a product probability measure
$\mu$ on $\{0,1\}^E$ (with $E$ a finite set)
and decreasing events $\A_i\sub \{0,1\}^E$, one has
$\mu(\cap \A_i)\geq \prod \mu(\A_i)$).
Thus (except for a factor $(1-o(1))$)
the probability that $\GG$ contains exactly $M$ triangles and no two of these
share a vertex is at least the product of the
right sides of \eqref{Scount} and \eqref{PrRS}, which is easily seen to be
$\exp[-O(m^3p^3)]$.
\qed

\section{Remarks}\label{Remarks}

For a fixed graph $H$ write $\xi=\xi_H(n,p)$ for the number of
(unlabeled, say) copies of $H$ in $G(n,p)$
and set $\mu=\mu_H(n,p)=\E \xi$.
A beautiful result of Janson, Oleszkiewicz and Ruci\'nski \cite{JOR}
says
\beq{jor}
\Pr(\xi > (1+\eta)\mu) < \exp[-\Omega_{H,\eta}(M(H,n,p))],
\enq
where we omit the definition of the parameter $M$
(understanding of which is one of the
main concerns of \cite{JOR}), but
just mention that
(i) if not zero, the probability in \eqref{jor} is bounded below by
$\exp_p[O_{H,\eta}(M(H,n,p))]$, and (ii)
except in the uninteresting case $p< n^{-2/(r-1)}$,
$M(K_r,n,p) $ is $\Theta(n^2p^{r-1})$.
In particular \eqref{jor} includes the result of \cite{Kim2} mentioned in Section \ref{Intro}
(and a bit more, since \cite{Kim2} requires $p>n^{-1}\log n$).

Of course \eqref{jor} is sharp when $H=K_2$, but we guess, perhaps
optimistically,
that this is the only case where the analogue of Theorem \ref{T1}
fails; more precisely, for any connected $H\neq K_2$ and
$p $ not too small, we should have
\beq{guess}
\Pr(\xi(H,n,p)>(1+\eta)\mu_H(n,p)) < p^{\Omega_{H,\eta}(M(H,n,p))}.
\enq
At this writing we think we can at least push the present
argument to
prove \eqref{guess}
for complete
graphs.
We won't try to say here what the lower bound on $p$ should be in general.
For $K_r$ it should be
$
n^{-2/(r-1)}(\log n)^{2/[(r-1)(r-2)]}
$,
which is (essentially) where the lower bound in (i) above
becomes larger than the bound in \eqref{guess}.

\bn
{\bf Added in proof.}
We have now completely settled the problem for cliques
and believe we know what should happen for general $H$ and $p$;
these items will appear in \cite{DK2}.

\bn
{\bf Acknowledgment.}
We would like to thank one of the referees for an exceptionally
careful reading.

\bn
Department of Mathematics\\
Rutgers University\\
Piscataway NJ 08854\\
rdemarco@math.rutgers.edu\\
jkahn@math.rutgers.edu

\end{document}